\documentclass[11pt,a4paper]{article}

\usepackage{amsmath}  
\usepackage{amssymb}  
\usepackage{latexsym} 
\usepackage{units}
\usepackage{hyperref}
\usepackage{tensor}
\usepackage{multirow}
\usepackage{setspace}             
\usepackage{fancyhdr}            
\usepackage{datetime}
\usepackage{array}
\usepackage{lscape}
\usepackage{rotating}
\usepackage[english]{babel}
\usepackage{txfonts,epsfig,graphicx,float}
\usepackage{epstopdf}
\usepackage[font=small,labelfont=bf]{caption}
\usepackage{subfigure}
\usepackage{enumerate}
\usepackage[table]{xcolor}
\usepackage{booktabs}

\usepackage{titlesec}

\titleformat*{\section}{\large\bfseries}
\titleformat*{\subsection}{\large\bfseries}
\titleformat*{\subsubsection}{\large\bfseries}
\titleformat*{\paragraph}{\large\bfseries}
\titleformat*{\subparagraph}{\large\bfseries}

\begin{document}

\begin{center}
{ \LARGE Sustainable theory of a logistic model - Fisher Information approach\\[0.3cm]}
\bigskip
\bigskip
\begin{flushleft}
{\bf{Avan Al-Saffar and Eun-jin Kim}}\\
{\small {School of Mathematics and Statistics, The University of Sheffield, Sheffield, South Yorkshire S3 7RH, UK.}}
{ \small {E-mail: asal-saffar1@sheffield.ac.uk and e.kim@sheffield.ac.uk}}
\end{flushleft}
\bigskip
\end{center}
 {\bf {Abstract}.} Information theory provides a useful tool to understand the evolution of complex nonlinear systems and their sustainability.  In particular, Fisher Information (FI) has been evoked as a useful measure of sustainability and the variability of dynamical systems including self-organising systems. By utilising FI, we investigate the sustainability of the logistic model for different perturbations in the positive and/or negative feedback. Specifically, we consider different oscillatory  modulations in the parameters for positive and negative feedbacks and investigate their effect on the evolution of the system and Probability Density Functions (PDFs). Depending on the relative time scale of the perturbation to the response time of the system (the linear growth rate), we demonstrate the maintenance of the initial condition for a long time, manifested by a broad bimodal PDF. We present the analysis of FI in different cases and elucidate its implications for the sustainability of  population dynamics. We also show that a purely oscillatory growth rate can lead to a finite amplitude solution while self-organisation of these systems can break down with an exponentially growing solution due to the fluctuation in negative feedback. \\

\section{Introduction}
\label{back}
 ~~~~~~~Nonlinear dynamical systems have been widely used as simple models for complex phenomena, e.g. in environmental, astrophysical and geophysical, and biological systems. In particular, the utility of such models in understanding biosystems has grown significantly in recent years as ever-improved experimental data has become available. A logistic model, first proposed by Verhults  to describe the growth of human populations in 1838 \cite{Sinha94, Loyd13}, is  one of the most popular models for the growth in biologist systems, e.g. bacteria, tumor cells, etc (\cite{Weirong05} and references therein). As a mean field equation, the logistic model describes the time-evolution of macroscopic  (large-scale) variables where the overall effect of micro-scale (small-scale) variables is incorporated by control parameters for the positive and negative feedback.  The merit of this model lies in the simplicity in the incorporation of the two complementary effects of a positive feedback which drives the growth  and a negative feedback which regulates its growth, thereby serving the simplest model for a self-regulated system where the growth is regulated within a system. The balance between the positive
and negative feedback leads to a stable equilibrium point (the so-called carrying capacity), to which a system approaches in a long time limit, regardless of the initial condition. Thus, a unique value of a carrying capacity can be viewed as a loss of the memory of the initial points. \\

Many researchers have extended the logistic model to include perturbations in the model parameters for feedback by periodic or random modulation (e.g. \cite{WANG11,AI03,WANG13,Kim16, Lu12,Liu12,Kim15,Prigogine84}) or to couple to the evolution of other systems (e.g. \cite{WANG09}). In particular, the possibility of bimodal Probability Density Function (PDF) was demonstrated in the presence of a correlation between a multiplicative noise (for the growth rate) and additive noise. The purpose of this paper to revisit this logistic model in view of the sustainability for different perturbations.  We compute PDFs for  different modulation in the model parameters and elucidate fundamental mechanisms determining the shape of PDFs. In particular, we demonstrate that when the characteristic time scale associated with the perturbation is much shorter than the system's response time, the system  maintains a long-term memory of initial conditions, thereby leading to a broad bimodal distribution.  The sustainability of a system in different cases is examined by computing FI. 
To test the stability of the most sustainable state inferred from our analysis of the FI, we add an additive noise to our system and test the resilience of our system to the environmental perturbation, modeled by the additive noise.\\

We note that the effect of fluctuating parameters have been studied in other dynamical systems \cite{Mohamed14a,Priede10,Cabezas02,Cabezas03, Mohamed14b,Shastri06,Mayer06,Sood13} while the dynamics of such systems has hardly been investigated from the perspective of information theory. Simplicity of the logistic equation enables us to undertake a systematic investigation in this regards.
The remainder of the paper is organised as follows. We introduce our model in \S 2 and present PDFs and Fisher information in \S 3 and \S 4, respectively, when the model parameters for both positive and negative feedback have the same fluctuations. In \S 5, we test the stability of our systems by adding an additive force.  Section 6 summarises the results for different types of modulation of the model parameters. Conclusions are provided in \S 7. \\

\section{Model and Motivation}
~~~~~~~ We  consider a population $x$  $(>0)$ and its logistic equation in the following form:
\begin{eqnarray}
  \frac{dx}{dt} =  N x \left ( 1- \frac{x}{K} \right ). \quad \quad
  \label{eq1}   
\end{eqnarray}
Here, $ N $ is the net growth rate, and $ K $ $(>0)$ is the carrying capacity of the system representing the maximum population size that can be supported by the system. The linear term  $Nx$ with $N>0$ represents a positive feedback while the nonlinear term  $N x^{2}/K$ represents a negative feedback. 
We note that regardless of the initial value of $x(t=0)= x_{0} $, $x$ reaches the carrying capacity $K$ as $t \to \infty$ for a constant $N>0$. \\

Compared to the case when the linear growth rate is constant or contains fluctuations in the absence/presence of an additive noise, it is less well understood what happens when the model parameter for the negative feedback contains fluctuations. Fluctuations in negative feedback can provide an interesting mathematical model for the loss of self-regulation, e.g. in biosystems (e.g. \cite{Kim16,Kim15}). While we comment on the cases where the model parameter for only positive or negative feedback contains fluctuations  in \S 6, of our particular interest in this paper is the case where the perturbation in positive and negative feedback is strongly correlated. Specifically, in  \S 2-4, we focus on the case of the following periodic modulation:
\begin{eqnarray}
  N = B+ N_{0} \sin (\omega t),
\label{eq2}
\end{eqnarray}
where $ B $ is a constant growth while $ N_{0}$ and  $\omega$  are the amplitude and frequency of the modulation. Since  $t$ and $x$ can always be rescaled by $N_{0}$ and $K$, respectively, we fix the value of $N_{0}$ and $K$ to be $N_{0}=5$ and $K=10$ and further focus on the case $B=0$ to study the effect of $\omega$ on the response of the logistic system. As shall be shown shortly, one of the consequences of the same  fluctuations in positive and negative feedbacks is the maintenance of an initial condition and bimodal distribution.\\
 The exact solution to Eq. (\ref{eq1}) with Eq. (\ref{eq2})  is easily found as:
\begin{eqnarray}
 x(t) = \frac{-K x_{0} \exp{\biggl(B t + \frac{N_{0}}{\omega} (1 - cos(\omega t))\biggr)}}{(x_{0} - K) - x_{0} \exp{\biggl(B t + \frac{N_{0}}{\omega} (1 - cos(\omega t))\biggr)}},
 \label{eq3}
\end{eqnarray}
where $x_{0}$ is the initial value of $x$ at $t=0$. The case $B=0$ represents the case where the linear growth rate has a zero average and fluctuates between $|N_{0}|$ and $-|N_{0}|$ in time, and is an interesting model for systems where a growth is strongly inhibited as in the case of bacteria under the action of antibiotics, etc. Thus,  we take $B=0$ and show the typical time history of $x(t)$ in Fig. 1 for different values of $\omega$ and $x_{0}$.
\begin{figure}[ht]
 \centering
           \includegraphics[width=0.9\textwidth]{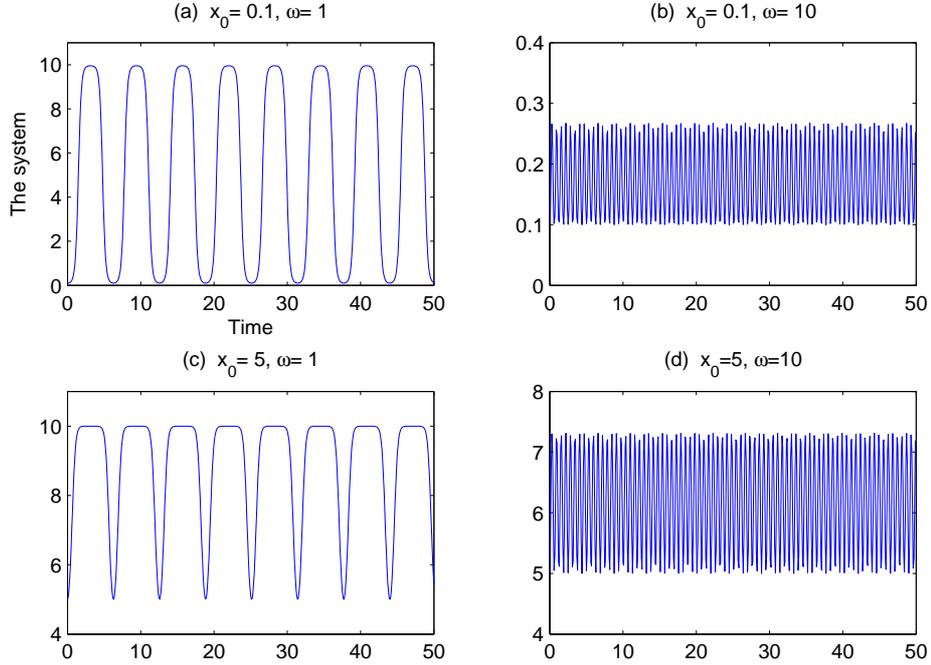}
           \label{Differentinitialconditions}
\captionsetup{justification=centering}
\caption{Time trace of $x(t)$ for different values of $x_{0}= 0.1, 5$ and $\omega=1,10$. For a small value of $\omega$, $x(t)$ tends to reach the carrying capacity $ K= 10$ while for large $\omega$, $x(t)$ maintains the initial condition. \\[0.4cm]}
\end{figure}

For a sufficiently small $\omega$, the time-scale of the perturbation becomes much larger than the system's response time (i.e. the mean square root  value of the growth rate), permitting enough time for $x$ to reach the carrying capacity, regardless of $x_{0}$. In comparison, for sufficiently large $\omega$ such that the perturbation occurs on time scales much shorter than the growth time, $x$ starting far from $x=10$ can never reach $x=10$ due to 
frequent periodic change in $N$, staying near $x=x_{0}$. To  demonstrate this cross-over between the case $ x \rightarrow K $ and $ x \rightarrow x_{0}$ for large $t$ in detail, we show the maximum of value of $x$ (in time) for different values of different $ \omega $ in Fig. 2, where the $ x$ and $y$ axes represent  $ \omega $ and the maximum values of $x$, respectively. Panels (a) and (b) show the maximum values of $x$ when $x_{0}= 0.1$ and $x_{0}= 5$, respectively.  To highlight the detailed feature for small $\omega$, the same figures in panels (a)-(b) are shown in log-log scale in panels (c)-(d), respectively.
From this, we observe a general tendency of the maximum $x$ monotonically decreasing as $\omega$ increases.
\vspace{10cm}
\begin{figure}[ht]
 \centering
           \includegraphics[width=0.9\textwidth]{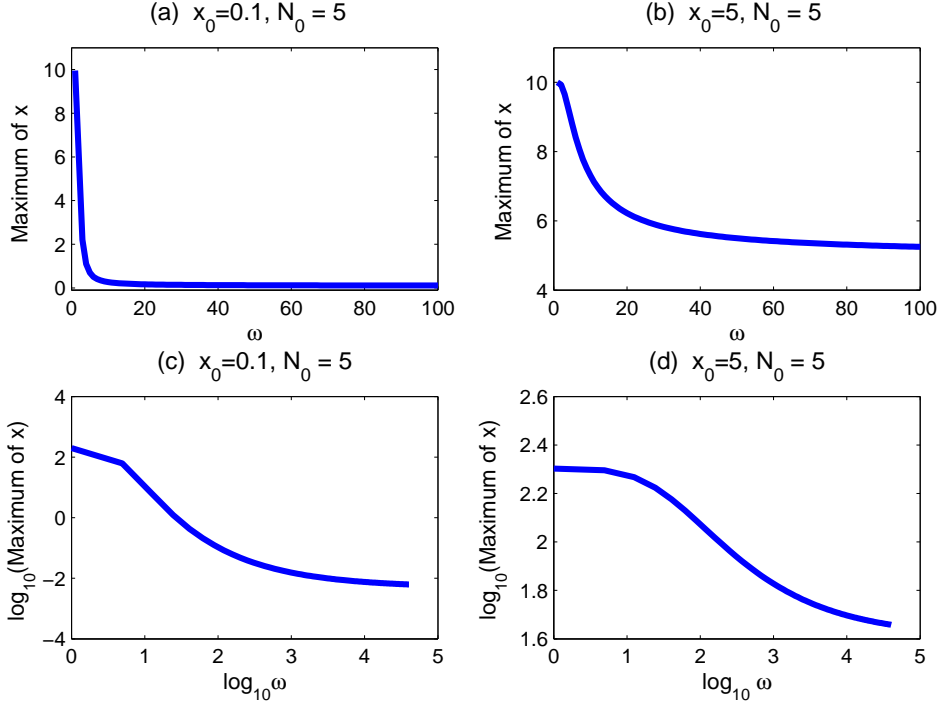}
           \label{DifferentinitialconditionsFixN0varyomega}
\captionsetup{justification=centering}
\caption{ The maximum value of $ x $ as a function of $\omega$ for $N_{0}= 5$ and $K=10$ by using different initial conditions $x_{0}= 0.1$ in panels (a) and (c) and $x_{0}= 5 $ in panels (b) and (d). (c)-(d) are shown in log-log scales .\\[0.4cm]}
\end{figure}\\
The main difference between the two cases with $x_{0}= 0.1$ and $x_{0}=5$ is a much steeper decrease in the maximum $x(t)$  for $x_{0}=0.1$ than for $x_{0}=5$. As the maximum $x(t)$ is obtained by the approach to the carrying capacity, the steep drop in the maximum $x(t)$ represents the inability of the system to reach this carrying capacity when the control parameter changes too rapidly in time. In this case, $x(t)$ does not deviate far from its initial value, effectively, leading to the maintenance of the memory of its initial value. This is consistent with the results shown in  Fig. 1.

\section{ Probability Density Function}

~~~~~~~~ We now examine the effect of $\omega$ and $x_{0}$  on Probability density functions (PDFs). To this
end, we compute PDF  by relating  the probability of observing the system at a particular value of $x$  to the amount of time the system state spends at $x$ \cite{KAMTHAN99} through conservation of the probability:
\begin{eqnarray}
       p[x] \;dx = p[t]\;dt. 
       \label{eq4}
\end{eqnarray}
Since  $t$ is a continuous variable with a uniform probability density: 
\begin{equation}
p[t] = constant = A, 
\label{eq5}
\end{equation}
we can obtain PDF of $x$ from Eqs. (\ref{eq4})-(\ref{eq5}) as:
\begin{eqnarray}
       p[x] = p[t]\; \biggl|\frac{dt}{dx}\biggr| = A\;  \biggl|\frac{dt}{dx}\biggr| = \frac {A}{u},
       \label{eq6}
\end{eqnarray}
where
\begin{equation}
       u = \frac{dx}{dt}.
       \label{eq7}
\end{equation}
The PDFs are shown for different values of $\omega$ in Fig. 3.
\begin{figure}[ht]
 \centering
           \includegraphics[width=0.9\textwidth]{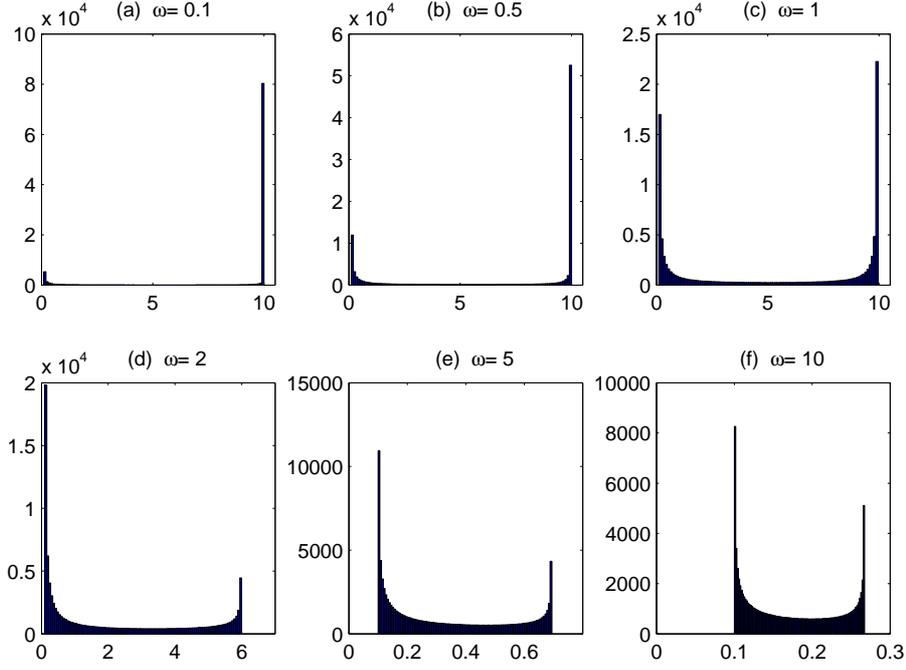}
           \label{Histogram}
\captionsetup{justification=centering}
\caption{PDF of $x(t)$ for  $N_{0}= 5$, $K=10$, $x_{0}= 0.1$.  Different values of $\omega=0.1,0.5,1,2,5,10$ are used in panels (a)-(f). A bimodel PDF is observed for all the cases.\\[0.4cm]}
\end{figure}\\

In Fig. 3, we observe a bimodal PDF for all the cases, with different distance between the two peaks. This bimodal distribution results from the maintenance of the initial condition $x_{0}=0.1$  against the tendency of $x$ approaching a carrying capacity ($10=K$), as noted previously. Specifically, for small $\omega \ll N_0/2\pi$ where the time-scale of the perturbation is much larger than the growth time $1/N$, $x$ reaches the carrying capacity, regardless of $x_{0}$, leading to the two peaks at $x=0.1$ (initial condition) and $x=10$ ($=K$ the carrying capacity).  In comparison, for sufficiently large $\omega \gg N_0/2\pi$ such that the perturbation occurs on time scales much shorter than the growth time (in root mean square value), $x$ starting far from $x=10$ can never reach $x=10$ due to frequent periodic change in $N$, leading to the formation of a very narrow distribution near $x=x_{0}$. This narrow PDF near $x_{0}$ manifests the maintenance of the initial condition when the perturbation occurs much faster than the system's response time. Between these two extreme cases, the bimodal PDF with the  largest distance between the two PDF peaks appears for the parameter $N_{0} / \omega=5$. It is interesting to observe the gradual shift of the population from the right PDF peak to the left PDF peak with the increase in $\omega$, followed by the narrowing of the PDF. That is, the narrowing of the PDF occurs after the left PDF peak around $x_{0}=0.1$ has grown taller than the right PDF peak. \\
\vspace{10cm}
\begin{figure}[ht]
 \centering
           \includegraphics[width=0.9\textwidth]{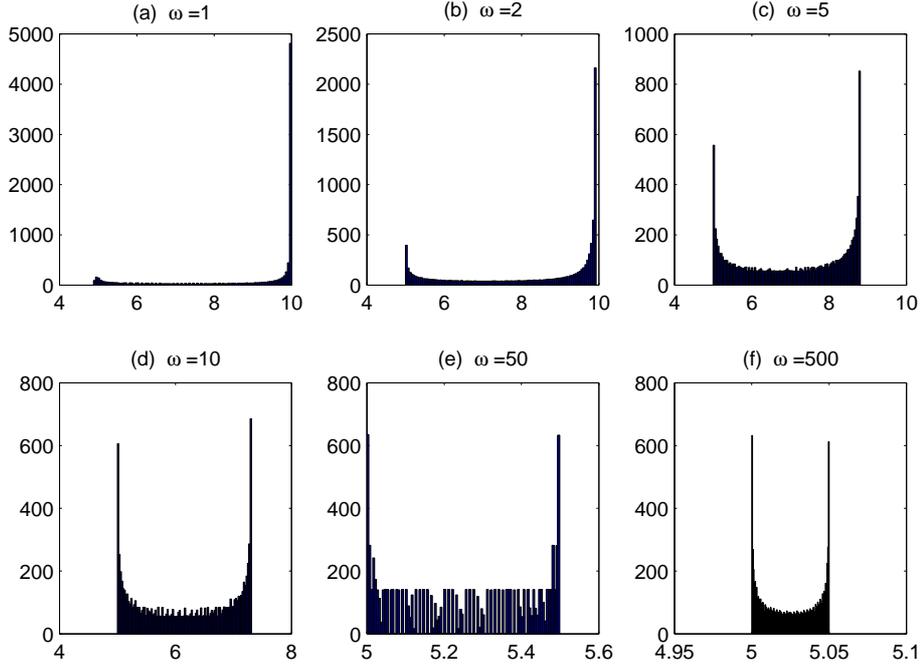}
           \label{Histogramx05}
\captionsetup{justification=centering}
\caption{ PDF of $x(t)$ for $N_{0}= 5$ and $ x_{0}= 5 $ by using different values of $\omega$ in panels (a)-(f). \\[0.4cm]}
\end{figure}\\

To demonstrate how PDF depends on $x_{0}$, we show another case in Fig. 4 by using the initial value  $x_{0}=5$ much closer to the carrying capacity. Similarly to the case $x_{0}=0.1$ in Fig. 3, Fig. 4 for $x_{0}= 5$  demonstrates a bimodal PDF for all cases; for small $\omega$,  $x$ reaches the carrying capacity while for large $\omega$, $x$ starting far from $ x=10$, can never reach $x=10$ and only fluctuates around $x_{0}= 5$, and the distance between the two peaks becomes shorter as $\omega$ increases. However, in contrast to Figs. 3, there is no significant growth of left PDF peak around $x_{0}=5$ for any $\omega$ and the narrowing of the PDF in Fig.4  occurs when the right PDF peak is still larger than the left PDF peak.\\
That is, for the initial condition $x_{0}=0.1$ (much less than the carrying capacity), there is an optimal value of $\omega$, which can maintain the distinct bimodal PDF while for the initial condition $x_{0}=5$ (close to the carrying capacity), such an optimal value of $\omega$ does not exist. The implication of the existence of such optimal value of $\omega$ will later be related to the utility of FI as a measure of the sustainability.

  \section{Fisher information}
~~~~~~~ Results shown in previous sections highlight a significant change to the logistic model due to periodic modulation in model parameters. 
 In this section, we examine this effect from the point view of  FI. FI is a function of the variability (order) of the observations such that low variability (strong order) leads to high FI. That is, a PDF bias to particular $ x $ values has higher FI whereas high variability (low order) with a lack of predictability of values of $ x$ leads to small FI (e.g. ``unbiased'' PDF).  This is demonstrated in Fig. 5.  The previous work suggested the following sustainability hypothesis: {\emph{``sustainable systems do not lose or gain Fisher information over time}}'' \cite{Cabezas03,Frieden95,Rico-Ramirez10}.
\begin{figure}[H]
 \centering
           \includegraphics[width=0.7\textwidth]{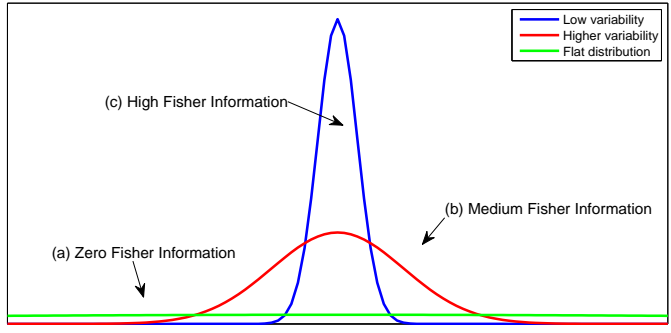}
           \label{Disorder}		
\vspace{-5pt}
\captionsetup{justification=centering}
\caption{(c) A steeply sloped PDF of $x$ with a large \\ FI (High gradients), (b) a PDF of $x$ with a smaller FI (a) a uniform PDF with zero FI.\\[0.4cm]}
\end{figure}
One of the utility of Fisher's information measure has been in the development of the basic theory of sustainability, for instance, in order to determine whether the system is sustainable or not \cite{Rico-Ramirez10,Suidan12} in diverse physical systems  (see \cite{Sanchez-Moreno11}  and references there in). We recall that FI is a very special uncertainty measure;  in contrast to a global measure of uncertainty (e.g., variance, or Shannon's entropy), FI strongly depends on the  gradient of PDF, consequently, is sensitive to the local oscillatory character of the PDF and relabeling \cite{Shastri06,Sanchez-Moreno11,Fath05}.\\

By following  Cabezac and Fath \cite{Cabezas03}, for a single variable  $x$, FI is calculated from the PDF of $x$, $p(x,t)$, as follows:\footnote{We note that Eq. (\ref{eq8}) can be extended to $n$-dimensional system. }
\begin{eqnarray}
        I = \int \frac{1}{p(x)}  \biggl(\frac {dp(x)}{\partial x}\biggr)^{2}\, d x.
        \label{eq8}
\end{eqnarray}
We compute the time averaged I by using Eqs. (\ref{eq6}), (\ref{eq7}) together with
\begin{eqnarray*}
       \frac {dp [x]}{dt} = - \frac {A}{u^{2}}\; \frac{du}{dt},
\end{eqnarray*}
in Eq. (\ref{eq8}) as follows:
\begin{eqnarray}
        I = \frac{1}{T} \int_0^T \frac{1}{(u(t))^{4}} \;  {\biggl({\frac{du}{dt}}\biggr)^{2}} \; d t = \int_0^T  \frac{1}{A}\, \biggl(\frac{dp(x)}{dt}\biggr)^{2}\; d t.   
\end{eqnarray}
Here, $I$ is the FI averaged over  the total time duration $T$; $A$  is a normalization constant. 
In the following, we investigate the sustainability/variability of our system by computing FI for different cases \cite{Cabezas03}. We use the same values  of $N_{0}=5$ and $K=10$, as before, and present FI for different values of  $\omega$ and for the two initial $x_{0}=0.1, 5$.\\

For each case with the fixed parameter/initial values ($\omega$ and $x_{0}$), we compute $I$  by varying the total time duration $T$, for instance, by using $t=[0, 10]$ with $T=10$, $t=[0,20]$ with  $T=20$, and so forth and present FI as a function of $T$.
Fig. 6 shows FI against $T$ for different $\omega = 0.1,0.5, 1,2,5,10$ in panels (a)-(f) for the fixed $x_{0}=0.1$,  corresponding to the case shown in Fig. 3. Specifically, we use $1000$ data points for each panel for $T= 10 n$  ($n=1,2,3, 1000$). In each panel, we observe that FI initially undergoes transient state and approaches an asymptotic value for a sufficiently large $T$. The higher asymptotic value of FI can be observed for  $\omega = 1$ while a small value is observed for  $\omega = 0.1 $.
We show how this asymptotic value of FI varies with $\omega$ in Fig. 7. A notable feature of Fig. 7 is the presence of a distinct maximum of FI around $\omega  \sim 1$, and this is related to the existence of the optimal $\omega$ which maintains the two peaks in the bimodal PDFs, discussed in relation to Fig. 3. \\
\begin{figure}[ht]
\centering
           \includegraphics[width=0.9\textwidth]{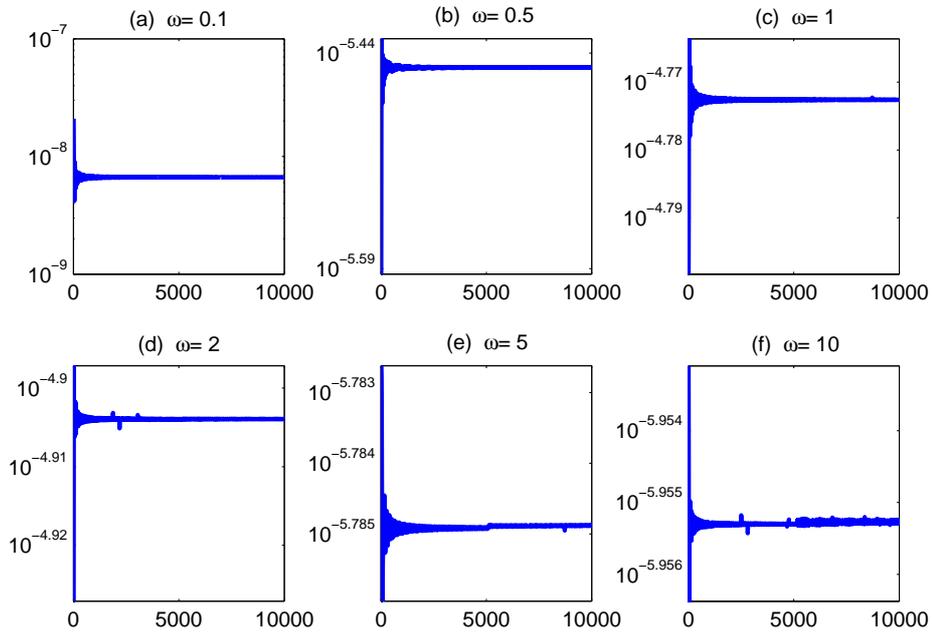}
           \label{Fisher information fixN0varyOmegasmallx0}	
\captionsetup{justification=centering}
\caption{FI against the total time $T$ for  $N_{0}= 5$, $ x_{0}= 0.1 $, $K= 10$. Panels (a)-(f) are different values of $\omega$.  We can observe the higher value of FI is when $\omega=1$. \\[0.3cm]}
\end{figure}
\vspace{20cm}
\begin{figure}
\centering
           \includegraphics[width=0.8\textwidth]{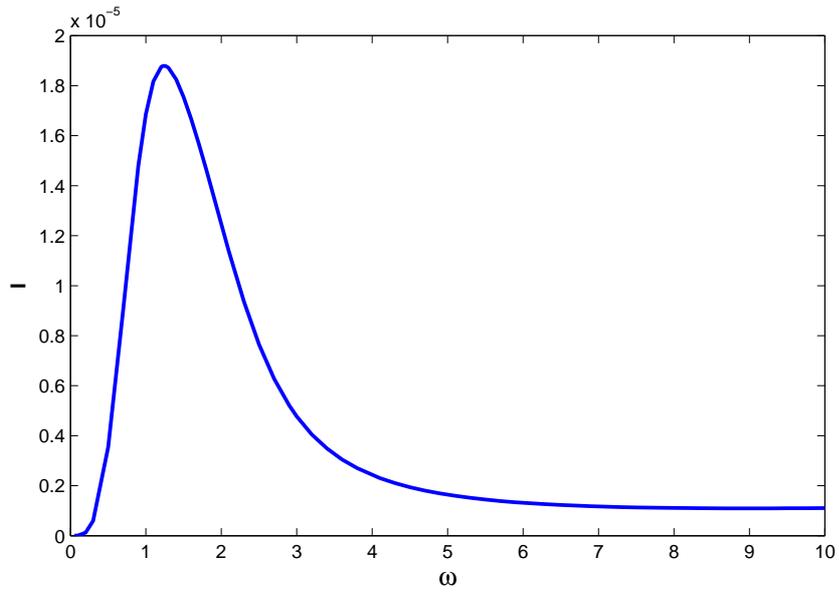}
           \label{Relation between omega and FI}	
\captionsetup{justification=centering}
\caption{Asymptotic value of FI against $\omega$ for $x_{0}= 0.1$. \\[0.3cm]}
\end{figure}\\

In the following, this distinct maximum in FI is shown to disappear in the case of $x_{0}=5$, the case corresponding to Fig. 4.
\vspace{15cm}
\begin{figure}[ht]
\centering
           \includegraphics[width=0.9\textwidth]{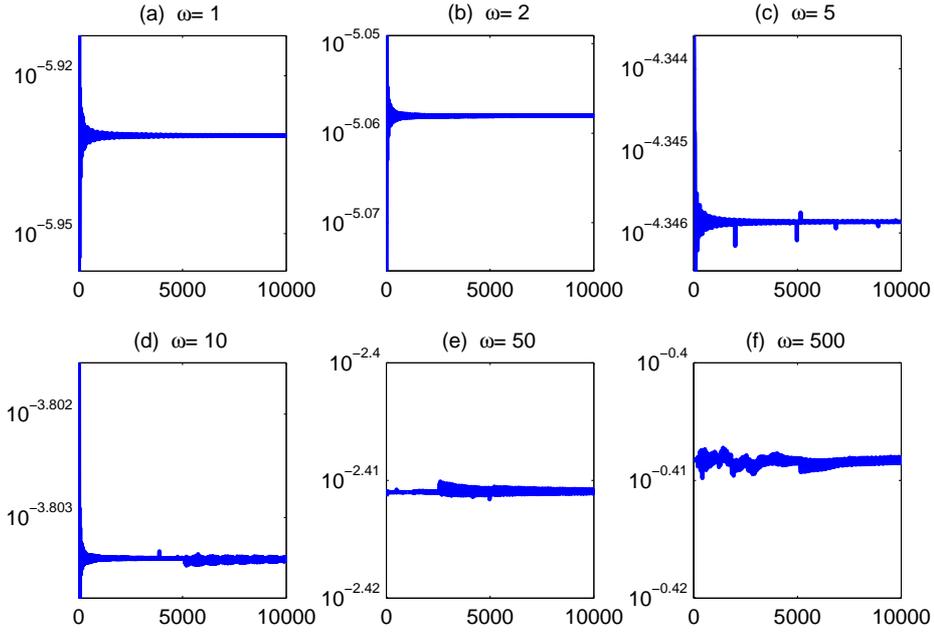}
           \label{Fisher information fixN0varyOmegasmallx0}	
\captionsetup{justification=centering}
\caption{Fisher Information fo $N_{0}= 5$, $K= 10$ and $ x_{0}= 5 $. Panels (a)-(f) are for different values of $\omega$. \\[0.4cm]}
\end{figure}

Figs. 8 and 9 show FI against $T$ for different values of $\omega$ and the asymptotic value of FI against $\omega$, respectively, for  $x_{0}=5$. Of notable is the monotonic increase of FI in Fig. 9, in a sharp contrast to Fig. 7. This represents that an optimal $\omega$ which maximise FI does not exist in this case; this is linked to the lack of the  two distinct peaks in bimodal PDFs, as discussed previously in relation to Fig. 4. \\

\vspace{15cm}
\begin{figure}[ht]
\centering
           \includegraphics[width=0.9\textwidth]{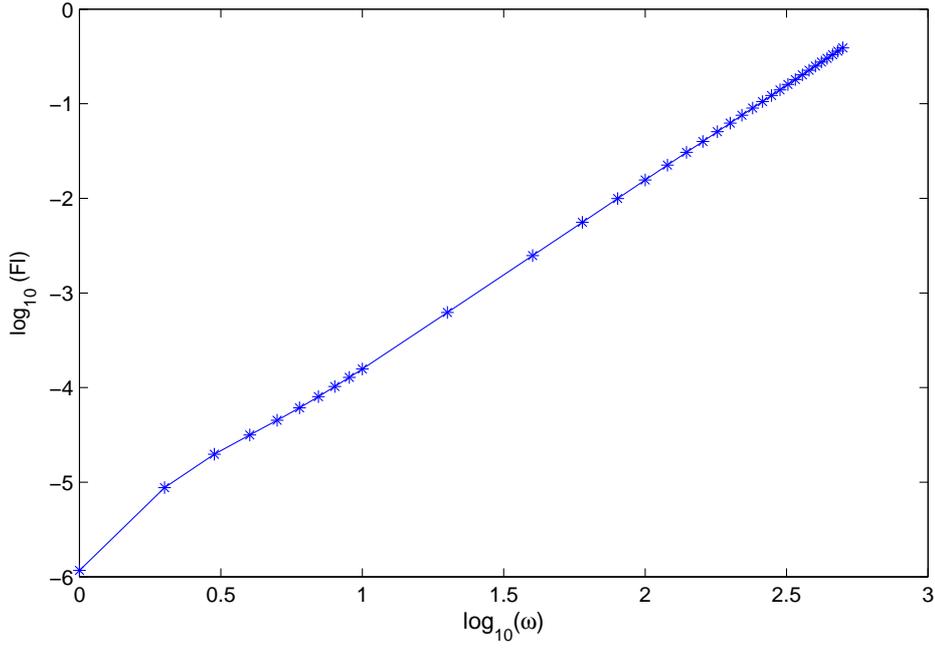}
           \label{Relation between omega and FI when $x_{0}= 5$}	
\captionsetup{justification=centering}
\caption{Asymptotic value of FI against $\omega$ for  $x_{0}= 5$. \\[0.4cm]}
\end{figure}
\section{Role of Fisher information as a measure of sustainability}
~~~~~~~  In previous sections we found  that the FI takes its maximum value around the optimal value  of $\omega \sim N_{0}/5 = 1$  when $x_{0}=0.1$ and $B=0$.  In order to test the sustainability of the optimal case with the maximum $I$, we examine the stability of  this optimal case by  adding an additive noise $B_{1} sin(\omega_{1} t)$ as follows:
\begin{eqnarray}
  \frac{dx}{dt} =  ( B + N_{0}  \sin(\omega  t) )  x   \biggl( 1- \frac{x}{K}\biggr ) +  B_1 sin(\omega_{1} t),
              \label{eq10}
\end{eqnarray}
and compare results with those obtained in non-optimal cases (e.g. $\omega=10$). 
We have explored different values of $\omega$, $B_{1}$ and $\omega_{1}$ and in the following, present the results for 
$\omega = 1$ (optimal case), $\omega=10$ (non-optimal case), $B_{1}=1,10$, $\omega_{1}=1,\sqrt{2}$ as example.
\begin{figure}[ht]
 \centering
           \includegraphics[width=0.9\textwidth]{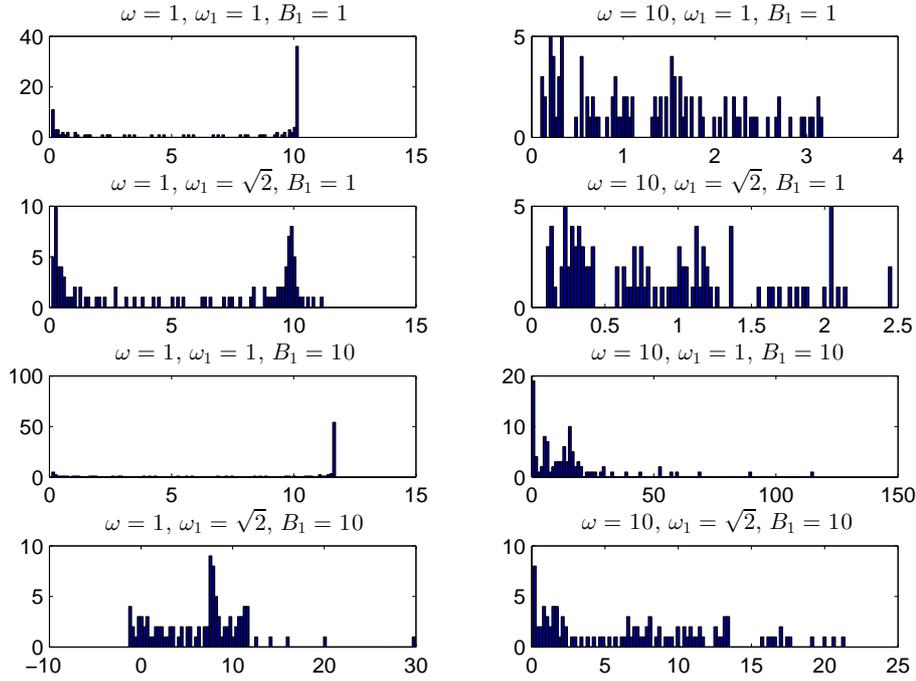}
           \label{FIcase-PDFadditivenoise}	
\vspace{-15pt}
\captionsetup{justification=centering}
\caption{PDF of $x$ for $x_{0}= 0.1$ and $N_{0}= 5$, with additive noise. Left panels are for the optimal case ($\omega=1$) while
right panels are for the non-optimal case ($\omega=10$). The PDFs with optimal  value $N_{0} =5 \omega$ in left panels  are more resilient to the additive noise than the PDFs for $\omega=10$.  \\[0.4cm]}
\end{figure}\ \
\vspace{10cm}

First, in Fig. 10, we show how PDFs are affected by different additive noise for $\omega=1$ in left panels and $\omega=10$ in right panels. In comparison with the PDFs in Fig. 3(c) and 3(f), respectively, we see that the overall change in PDFs is much less in the optimal case ($\omega=1$), suggesting that the optimal case with a large FI is less affected by the additive noise than in the non-optimal case. To strengthen this argument, we utilise the mean value as another measure to quantify the change of the system. Specifically, we compute the mean value without the additive noise (shown in Table 1) and the mean value after adding the additive noise, and quantify the ratio of the change in the mean value  as follows:\\
\begin{eqnarray}
  \frac{ \textbf{Mean value without additive noise} -  \textbf{Mean value with additive noise}}{ \textbf{Mean value without additive noise}} x 100\%.
\end{eqnarray}\\

The results are shown in Table 2. 
 We can see that in the optimal case $\omega = 1$, the ratio of change in the mean value for different additive noises is much less than that in the non-optimal case ($\omega = 10$). \\

To complete our investigate on the implication of FI for sustainability, we have also performed similar experiments for $x_{0}= 5$ by adding an additive noise of different amplitude and frequencies, and have found no obvious link between the value of FI and sustainability. 
This is due to the lack of maximum in FI for this initial condition. From these, we propose that the FI is a useful measure in the case when the FI has a distinct maximum (related to the presence of the two distinct bimodal PDF peaks).\\

\begin{table}[ht]
\centering
\begin{tabular}{|c|c|} \hline  
Mean value in  optimal case  & Mean value in  non-optimal case \ \\\hline  
5.5433 & 0.1733  \\\hline  
\end{tabular}
\caption{Mean value in the optimal and non-optimal cases without additive noise; $x_{0}=0.1$}
\label{rapidtable} \vspace{1cm} 
\end{table} 

\begin{table}[h]
\begin{center}
\begin{tabular}{c|c|c|c|l|}
\cline{2-5}
&  \multicolumn{2}{c|}{ Optimal case } 
& \multicolumn{2}{c|}{Non-optimal case} \\ \cline{2-5}
&  Mean value & \% Change (11) & Mean value & \% Change (11) \\ \cline{1-5}
\multicolumn{1}{ |c| }{$B_{1}=1, \omega_{1}=1$ } & 6.3927 & 15.3 \% & 1.3103 & 656.1 \%      \\ \cline{1-5}
\multicolumn{1}{ |c| }{$B_{1}=1, \omega_{1}=\sqrt{2}$} &  5.5515 & 0.2 \% & 0.8814 & 408.6 \%      \\ \cline{1-5}
\multicolumn{1}{ |c| }{$ B_{1}=10, \omega_{1}=1$} & 8.8326 & 59.3 \% & 15.2631 & 8707.3 \% \\ \cline{1-5}
\multicolumn{1}{ |c| }{$B_{1}=10, \omega_{1}=\sqrt{2}$} & 6.5509 & 18.2 \% & 7.2068 & 4058.6 \% \\ \cline{1-5}
\end{tabular}
\caption{\% change in mean value in the optimal and non-optimal cases with additive noise; $x_{0}=0.1$}
\label{rapidtable} \vspace{1cm} 
\end{center}
\end{table}
\section{Comments on different modulation}
~~~~~~~ We have so far focused on the case where the same modulation is applied to both positive and negative feedbacks. To complete our investigation, we now comment on the effect of the two different modulations. \\

\subsection{Case-1: Perturbation in the positive feedback}

~~~~~~~~We consider a periodic modulation in the parameter for the positive feedback for the constant model parameter in the negative feedback. Specifically, we consider:
\begin{equation}
  \frac{dx}{dt} =  [ B + N_{0}  \sin(\omega t )] x  - \frac {C x^{2}}{K},
           \label{eq11}
\end{equation}
where the values of $B$, $C$, and $K$ are kept constant.  In Fig. 11, we illustrate
the effect of different values of  $\omega$ and $N_{0}$ on PDFs for $ B = 0$, $K= 10$, $C= 1$, and $x_{0}= 0.1$. By taking $B=0$, we are again modelling
the case where the growth is strongly inhibited and is driven only by fluctuations. Even when the linear growth rate has zero average, we observe the excitation of the finite amplitude solution, similar to the result in \cite{Kim15}. This finite amplitude solution leads to PDFs centered around the initial position $x_{0}=0.1$ with a single peak, as shown in Fig. 11. That is, in contrast to a bimodal PDFs in the previous sections, we observe a unimodal PDF in all cases. This reflects the main effect of a multiplicative noise in driving a unimodal PDF. The  width of PDFs near $x=x_{0}$ becomes narrower as $\omega$ increases, similarly to the behaviour of the bimodal PDFs in the previous sections. \\
\begin{figure}[ht]
 \centering
           \includegraphics[width=0.9\textwidth]{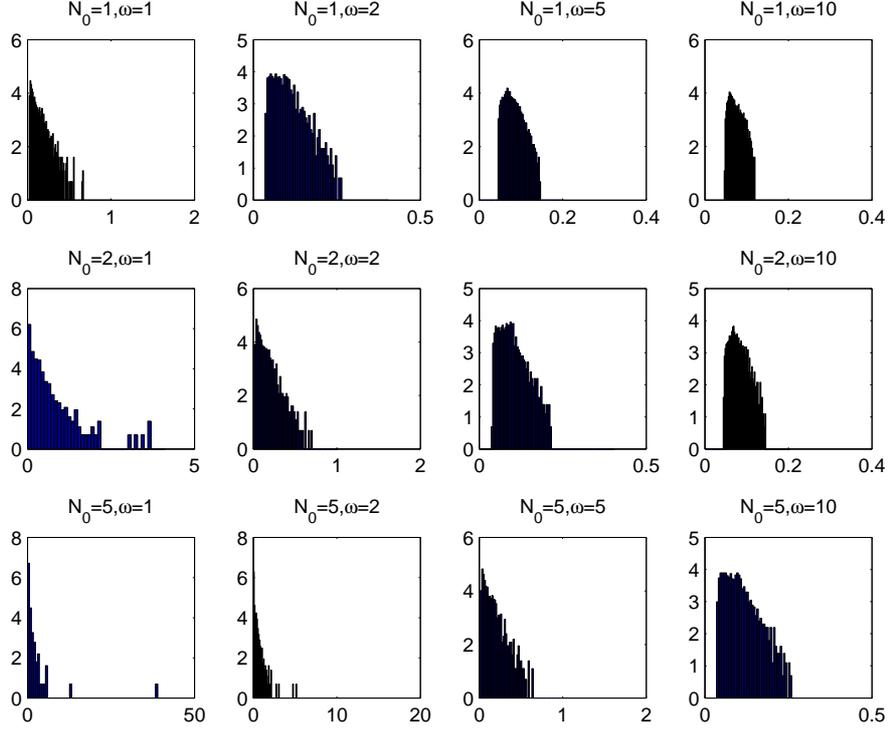}
           \label{Histogram2}
\captionsetup{justification=centering}
\caption{ PDFs of $x$ for Case-1 for different $N_{0}$ and $\omega$. $ x_{0}= 0.1 $,$ B = 0$, $K= 10$, $C= 1$. We observe notable  unimodal PDFs for all the cases.\\[0.4cm]}
\end{figure}\\

\subsection{Case-2: Perturbation in the negative feedback}

~~~~~~~~~ We finally consider that case where  fluctuations in the parameter are included only in the negative feedback as follows:
\begin{eqnarray}
  \frac{dx}{dt} =  C  x - \frac {[B + N_{0} \sin ( \omega  t )]x^{2} }{K}.
           \label{eq12}
\end{eqnarray}
The analytical solution to Eq. (\ref{eq12}) can be found as:\\
\begin{equation}
   x = \frac{K a b C x_{0}}{ K C a + N_{0} C b x_{0} c + N_{0} C \omega x_{0} + B x_{0} a (b-1)}, 
           \label{eq13}
\end{equation}\\
where
$$   a = C^{2} + \omega^{2},$$
$$  b = \exp(C t), $$
$$  c = C \sin(\omega t) - \omega \cos(\omega t). $$\\ 
As the amplitude of $N_{0}$ relative to $B$ increases, the solution starts growing exponentially as the nonlinear damping becomes ineffective (e.g. see \cite{Kim16,Kim15}).
\begin{figure}[ht]
 \centering
           \includegraphics[width=0.9\textwidth]{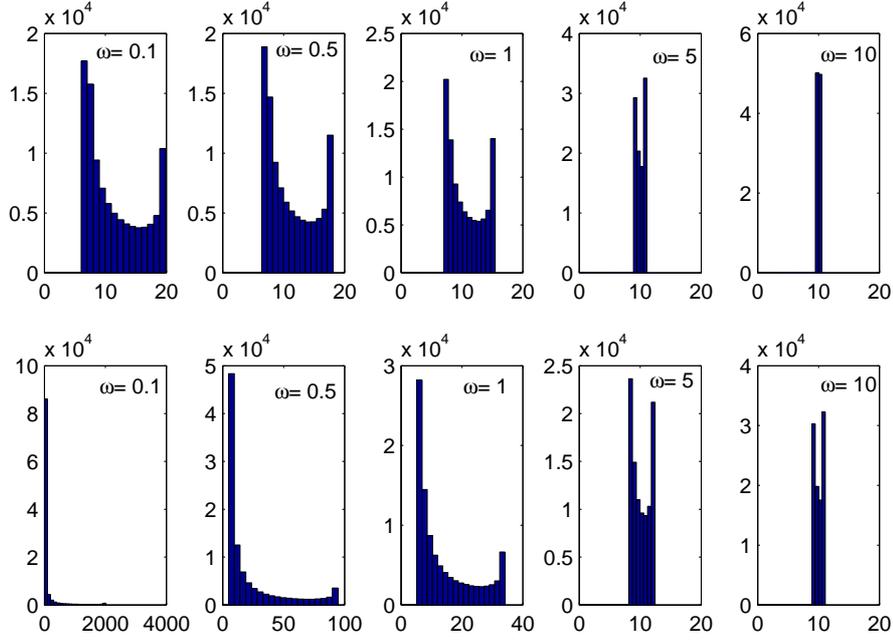}
           \label{Histogram3}
\captionsetup{justification=centering}
\caption{PDFs of $x$ for two different $N_{0}= 0.5$ and $1$ in the upper and lower panels, respectively. For all cases, $ x_{0}= 0.1 $, $B=1$, $K=10$, and $C= 1$, PDFs are unimodal in all cases.\\[0.4cm]}
\end{figure}\\ 
\vspace{5cm}

The resulting PDFs are shown in Fig. 12 for $N_{0}= 0.5$  in the upper panels and $N_{0}=1$ in the lower panels, respectively, for the same $ x_{0}= 0.1 $, $B=1$, $K=10$, and $C= 1$.  We observe that when $N_{0}= 1$, PDF becomes broader as $\omega$ decreases. The broadening of PDFs is related to the strong intermittency of $x$, manifested by the high-amplitude peaks, as $\omega$ decreases and can be seen from the time trace in Fig.13. In particular, we note that the solution grows exponentially for sufficiently large $N_{0}$ and
small $\omega$, as shown in panel (c) for  $\omega=0.1$ and $N_{0}=10$.
\begin{figure}[ht]
 \centering
           \includegraphics[width=0.8\textwidth]{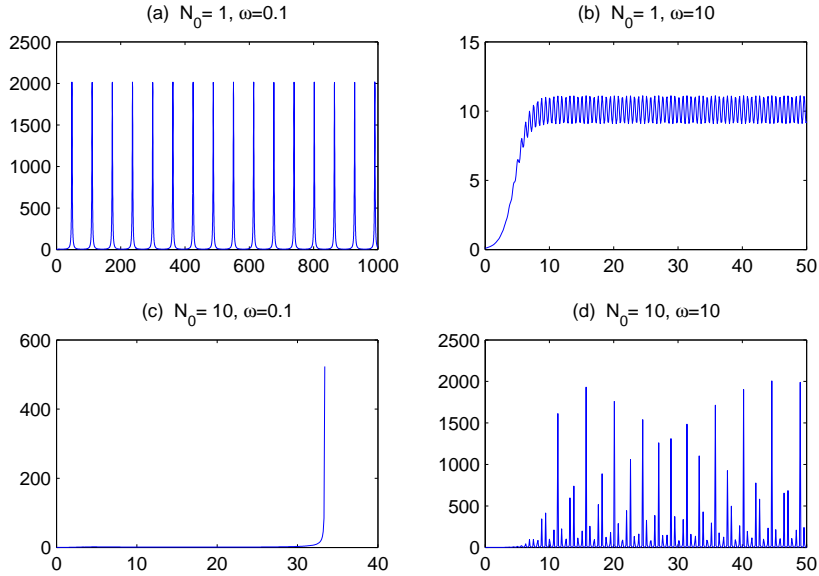}
           \label{Histogramcase-2}
\captionsetup{justification=centering}
\caption{$x(t)$ against $t$ for  $N_{0}=1,10$ and $\omega=0.1,10$. For all cases
 $x_{0}= 0.1$, $K= 10$, $B= 1$, $C= 1$
  \\[0.4cm]}
\end{figure}\\
\newpage 
\section{Conclusions and discussion} 
We have revisited the logistic model in view of the sustainability for different perturbations in the model parameters for the positive and/or negative feedback and investigated the effect of different modulation  and initial conditions. In particular, we demonstrated the possibility of the maintenance of a long-term memory of initial conditions when the characteristic time scale associated with the disturbance is much shorter than the system's response time,  as well as bimodal distributions. In the case of the same periodic modulation of the model parameters for the positive and negative feedback, for the initial condition far from the carrying capacity ($x_{0}= 0.1$), we found a distinct maximum value of FI for an optimal value of parameters $N_{0} \sim 5 \omega$ due to a broad bimodal PDFs with two distinct peaks.
In contrast, for $x_{0}= 5$, FI was shown to monotonically increases with $\omega$, with no distinct maximum. 
The sustainability of a system under different perturbation is examined by computing   FI from  PDFs. 
In particular, we found that FI is a useful measure of sustainability in the case when it has a distinct maximum, as a consequence of the presence of the two distinct bimodal PDF peaks. Our results could have interesting implications for understanding  the origin of  bimodal PDFs (e.g. \cite{PERSISTER}). That is, population of small size (corresponding to small $x_{0}$ in our model) can maintain a broad bimodal PDF for an optimal perturbation frequency with the maximum FI and thus has the best survival likelihood. As the optimal perturbation frequency occurs when its time scale is of order of the time scale of the linear growth rate (in root mean square value), it may well be that the population with such growth rate would have the best fitness. It would be of interest to extend our work to other systems such as a coupled logistic equations and a Gompertzian equation and study their implications  in future publications.
\bigskip
\bibliographystyle{myunsrt}
\bibliography{references}
 
\end{document}